\documentclass[12pt , reqno, fleqn]{amsart}
\usepackage[margin=1in]{geometry}
\usepackage{amssymb}
\usepackage{amsthm}
\usepackage{amstext}
\usepackage{amsbsy}
\usepackage{amscd}
\usepackage{enumerate}
\usepackage{chngpage}
\usepackage{mathtools}
\usepackage{amsmath}
\usepackage{hyperref}
\usepackage{tikz}
\usepackage{stmaryrd}
\usepackage{color}

\newtheorem{thm}{Theorem}[section]
\newtheorem{alg}[thm]{Algorithm}
 \newtheorem{cor}[thm]{Corollary}
\newtheorem{prop}[thm]{Proposition}
\theoremstyle{definition}
\newtheorem{defn}[thm]{Definition}
\newtheorem{ques}[thm]{Question}

\newtheorem{example}[thm]{Example}

\newcommand{\Z}{\mathbb Z}

\newcommand{\Q}{\mathbb Q}
\newcommand{\F}{\mathbb F}

\newcommand{\Fq}{\mathbb{F}_q}

\DeclareMathOperator{\PGL}{PGL}
\newcommand{\mb}[1]{\mathbb{#1}}
\newcommand{\mc}[1]{\mathcal{#1}}

\begin{document}
\title{An algorithm to count the number of caps in $\mathbb{P}^3(\mb{F}_q)$}
\author{Kelly Isham}
\date{}
\address{Department of Mathematics\\ 214 McGregory Hall\\
	Colgate University\\
Hamilton, NY}
\email{kisham@colgate.edu}
\maketitle

\begin{abstract}
	An $n$-cap in $k$-dimensional projective space is a set of $n$ points so that no three lie on a line. In this note, we provide an algorithm to count the number of $n$-caps in $\mathbb{P}^3(\F_q)$, which follows from our recent paper \cite{ish_arcsp3}. We then give exact formulas for the number of $n$-caps when $n \le 7$. The formulas are polynomial in $q$ when $n \le 6$ and quasipolynomial in $q$ when $n = 7$. 
\end{abstract}

\section{Introduction}\label{intro_sec}

An $n$-arc in the projective plane is a set of $n$ points so that no three lie on a line. We can generalize this notion to higher-dimensional space in a couple different ways.

First, we define an $n$-arc in $k$-dimensional space to be a set of $n$ points so that no $k+1$ lie on a hyperplane. This object is studied across mathematics in areas including coding theory, finite geometry, algebraic geometry, and linear algebra. In particular, $n$-arcs in $\mb{P}^{k}(\F_q)$ are closely related to MDS codes, which are codes that achieve the Singleton bound. One important conjecture in coding theory is the MDS conjecture, which, rewritten in the language of arcs, states that if $k \le q$, the size of the largest arc in $\mb{P}^{k}(\F_q)$ is $q+2$ if $q$ is even and $k \in \{2, q-2\}$, and is $q+1$ otherwise. Researchers from various areas of mathematics have made substantial progress toward proving the MDS conjecture; see \cite{balllavrauw} for a survey of large arcs.

Second, we define an $n$-cap in $k$-dimensional projective space to be a set of $n$ points so that no three lie on a line.  Notice that the definition of $n$-arc and $n$-cap coincide when $k=2$. In this paper, we will consider $n$-caps. 

In 1936, Erd\H{o}s and Turán \cite{erdost} conjectured that sets of integers with positive density must have arbitrarily long arithmetic progressions. Famously, Szemerédi \cite{szemerdi} proved their conjecture in 1975 using combinatorial methods. The Cap Set problem is an analogue, which asks for the largest $n$ so that an $n$-cap exists in $\F_3^k$ (e.g. the largest set which does not contain a three-term arithmetic progression).  More generally, we can ask for the largest $n$ so that an $n$-cap exists in $\F_p^k$. Substantial progress on these problems has been made, see \cite{peluse} for a summary of known results and methods. The best known upper bounds are due to Ellenberg and Gijswijt \cite{eg} in a groundbreaking paper from 2017. 

Some results about maximal caps are also known in $\mb{P}^{k}(\F_q)$. It is known that the maximal size of a cap in $\mb{P}^3(\F_q)$ is $q^2+1$ unless $q =2$.  Further, in 1955, Barlotti \cite{barlotti} and Panella \cite{panella} independently showed that if $q$ is odd or $q=4$, there is a one-to-one correspondence between elliptic quadrics and maximal caps in $\mb{P}^3(\F_q)$. If $q =  2^r$ and $r$ is odd, there exist maximal caps which are not elliptic quadrics \cite{tits}. Much less is known when $k > 3$.  For a recent survey on caps in $\mb{P}^{k}(\F_q)$, see \cite{hirschfeldthas}.

In this paper, we consider the counting version of the Cap Set problem. 

\begin{ques}
	Let $c_{n,k}(q)$ denote the number of $n$-caps in $\mb{P}^{k}(\F_q)$. For fixed $k$ and $n$, what is $c_{n,k}(q)$ as a function of $q$?
\end{ques}

It may be difficult to give an explicit formula for $c_{n,k}(q)$, so we can ask to classify the function. We use the hierarchy: polynomial, quasipolynomial, and non-quasipolynomial. A function $f(q)  \in \Q[q]$ is \emph{quasipolynomial} if there exist finitely many polynomials $g_0(q), g_1(q),$ $\ldots, g_{N-1}(q) \in \Q[q]$ so that $f(q) = g_i(q)$ whenever $q \equiv i \pmod{N}$.

When $k=2$, caps are equivalent to arcs. In this case, $c_{n,2}(q)$ is known explicitly when $n \le 9$ \cite{glynn, iss, kklpw}. This function is polynomial in $q$ when $n \le 6$ and is quasipolynomial in $q$ when $n \in \{7,8,9\}$. In recent work joint with several authors, we show that $c_{10,2}(q)$ is non-quasipolynomial in $q$, though we are unable to give an explicit formula. Most likely these results together give a full classification of $c_{n,2}(q)$. In fact, in the upcoming joint paper, we conjecture that $c_{n,2}(q)$ will continue to be non-quasipolynomial for all $n \ge 11$. 

One way to count $n$-caps in $\mb{P}^2(\F_q)$ is to count the opposite -- namely sets of $n$ points so that at least one subset of three points is contained in a line. Glynn provides a formula \cite[Theorem 3.16]{glynn} that shows we can do much better than this naive inclusion-exclusion count by studying special combinatorial objects that contain several three-point lines. 
 
\section{Connection to arcs in $\mb{P}^3(\F_q)$}
Motivated by Glynn's ideas, in a recent paper \cite{ish_arcsp3} we provide a formula to count $n$-arcs in $\mathbb{P}^3(\F_q)$ in terms of a small number of combinatorial objects. In this paper, we will derive a formula to count $c_{n,3}(q)$ that follows from the formula to count $n$-arcs in $\mb{P}^3(\F_q)$. In order to state the algorithm, we first need several definitions from \cite{ish_arcsp3}.

While our main results hold over $\mb{P}^3(\F_q)$ -- the unique finite 3-dimensional projective space of order $q$ -- we can first consider caps in a more general space. 

\begin{defn}\cite{bb}
	A \emph{linear space} $(\mc{P}, \mc{L})$ is a pair of sets where $\mc{P}$ denotes a set of \emph{points} and $\mc{L}$ denotes a set of \emph{lines} satisfying the following properties:
	\begin{enumerate}
		\item Every line in $\mc{L}$ is a subset of $\mc{P}$.
		\item Any two distinct points belong to exactly one line in $\mc{L}$.
		\item Every line in $\mc{L}$ contains at least 2 points.
	\end{enumerate}
\end{defn}

Let $2^X$ denote the power set of the set $X$.
\begin{defn} \cite[Definition 1.8]{ish_arcsp3}
	A \emph{planar space} is a triple of sets $(\mc{P}, \mc{L}, \mc{H})$ where $\mc{P}$ is the set of points, $\mc{L}$ is the set of lines, and $\mc{H}$ is the set of planes such that 
	\begin{enumerate}
		\item $\mc{L}, \mc{H} \subseteq 2^{\mc{P}}$.
		\item $(\mc{P}, \mc{L})$ is a linear space.
		\item Any three distinct non-collinear points lie on a unique plane.
	\end{enumerate}
\end{defn}

Linear and planar spaces are very general. For example, affine, projective, and Euclidean spaces satisfy these axioms. In order to count arcs in $\mb{P}^3(\F_q)$, we define maps that embed the combinatorial sets of points, lines, and planes into $\mb{P}^3(\F_q)$.

A \emph{strong realization} of a planar space $f = (\mc{P}, \mc{L}, \mc{H})$ is an injective mapping $\sigma: \mc{P} \hookrightarrow \mb{P}^3(\F_q)$ so that
\begin{itemize}
	\item $\sigma(L)$ is contained in a line of $\mb{P}^3(\F_q)$ if and only if $L \subset \mc{L}$ and 
	\item $\sigma(H)$ is contained in a plane of $\mb{P}^3(\F_q)$ if and only if $H \subset \mc{H}$.
\end{itemize} 
The first condition states that $\sigma$ \emph{preserves lines} and the second states that $\sigma$ \emph{preserves planes}. We let $A_f(4,q)$ denote the number of strong realizations of $f$ to remain consistent with the notation of \cite{ish_arcsp3} where $\sigma$ maps $\mc{P}$ into $\mb{P}^{k-1}(\F_q)$ for $k =4$. 

Two planar spaces $f= (\mc{P}_1, \mc{L}_1,\mc{H}_1)$ and $g = (\mc{P}_2, \mc{L}_2, \mc{H}_2)$ are \emph{isomorphic} if there exists a bijective mapping $\rho: \mc{P}_1 \rightarrow \mc{P}_2$ that preserves lines and planes.

In a linear or planar space, every two points determine a line. In a planar space, every three non-collinear points determine a plane. Thus we refer to a line containing at least three points a \emph{full line} and a plane containing at least four points a \emph{full plane}. Each planar space is determined uniquely by its full lines and planes, so we drop the word full and ignore two-point lines and three-point planes when writing a planar space. We say a point has \emph{index} $(i,j)$ if it lies on exactly $i$ (full) planes and $j$ (full) lines.

\begin{example}\label{pspace_ex}
	Consider the planar space $f = (\mc{P}, \mc{L} , \mc{H})$ on the points $\{1,2,3,4\}$ such that all four points lie on a common plane and the first three points lie on a line. The planar space can simply be given by $\mc{P} = \{1,2,3,4\}$, $\mc{L} = \{\{1,2,3\}\}$, and $\mc{H} = \{\{1,2,3,4\}\}$. Implicitly the line $\{1,2\}$ is contained in $\mc{L}$ since every two points determine a line, but we suppress this information when writing the set $\mc{L}$.
\end{example}

Notice that a cap can be defined as a planar space which contains no full lines. In order to count caps in $\mb{P}^3(\F_q)$, we will count the number of strong realizations of planar spaces so that $\mc{L}$ contains at least one full line, then use inclusion-exclusion. The objects remaining in the summation are special combinatorial objects called hyperfigurations.

\begin{defn} \cite[Definition 1.12]{ish_arcsp3}
	\label{hyperfig_def}
	A \emph{hyperfiguration} is a planar space on $n$ points such that for every point $P$, the index $(i,j)$ of $P$ satisfies at least one of the following:
	\begin{enumerate}
		\item $i \ge 4$
		\item $j \ge 3$
		\item $(i,j) = (3,0)$.
	\end{enumerate}
\end{defn}

\begin{example}
	The planar space from Example \ref{pspace_ex} is not a hyperfiguration since the index of point 1 is $(1,1)$.
\end{example}

This definition is not intuitive. The proof of Theorem \ref{p3count} relies on a lemma \cite[Lemma 2.3]{ish_arcsp3} that allows us to work inductively by removing a point of index $(i,j)$ from a planar space $f$. The lemma assumes that $(i,j)$ does not satisfy any of the three conditions given in the definition of hyperfiguration.

\begin{thm}\label{p3count} \cite[Theorem 1.13]{ish_arcsp3}
	There exist polynomials $p(q)$ and $p_h(q)$ in $\Z[q]$ such that
	$$
	\#\{n\textrm{-arcs in } \mb{P}^3(\F_q)\} = p(q) + \sum_{h} p_h(q) A_h(4,q)
	$$
	where the summation runs over all isomorphism classes $h$ of hyperfigurations on at most $n$ points. Moreover, there is an algorithm that produces $p(q)$ and $p_h(q)$ for each isomorphism class $h$.
\end{thm}

While not explicitly stated in \cite{ish_arcsp3}, the algorithm is inductive and produces a similar expression for $A_f(4,q)$ for every planar space $f$ on at most $n$ points. We state this stronger fact as a corollary.

\begin{cor}\label{Af_cor}
	For any planar space $f$, there exist polynomials $r(q)$ and $r_h(q)$ in $\Z[q]$ such that
	$$
	A_{f}(4,q) = r(q) + \sum_{h} r_h(q) A_h(4,q)
	$$
	where the summation runs over all isomorphism classes $h$ of hyperfigurations on at most $n$ points. Moreover, there is an algorithm that produces $r(q)$ and $r_h(q)$ for each isomorphism class $h$.
\end{cor}

\subsection{Main Results}
We begin by introducing an algorithm to count $c_{n,3}(q)$ in terms of a summation over all planar space functions on at most $n$ points. Let $\mc{N}_n$ be the set of all planar spaces on at most $n$ points for which the largest size of a subset in $\mc{L}$ is two, i.e., the planar spaces which do not contain a full line. We can understand $c_{n,3}(q)$ by counting the strong realizations of each planar space $f \in \mc{N}_n$, then summing.  This gives the following algorithm.

\noindent \begin{alg} \label{caps_alg} 	\leavevmode
	\begin{enumerate}
	
		\item Use Corollary \ref{Af_cor} to compute $A_f(4,q)$ for each planar space $f$ on at most $n$ points.
		\item Output $c_{n,3}(q) = \sum_{f \in \mc{N}_n } A_f(4,q).$
	\end{enumerate}
\end{alg}

Applying Algorithm \ref{caps_alg} leads to the following main theorem.

\begin{thm}
	There exist polynomials $s(q)$ and $s_h(q)$ so that 
	$$
	c_{n,3}(q) = s(q) + \sum_{h} s_h(q) A_h(4,q) 
	$$
where the summation runs over all isomorphism classes $h$ of hyperfigurations on at most $n$ points. Moreover, Algorithm \ref{caps_alg} produces $s(q)$ and $s_h(q)$ for each isomorphism class $h$.
\end{thm}

It is natural to wonder what the degrees of the polynomials $s_h(q)$ are. In Propositions \ref{prop_6caps} and \ref{prop_7caps}, we will show that for $n \le 7$, the polynomials $s_h(q)$ in $c_{n,3}(q)$ are constants when $h$ is a hyperfiguration on exactly $n$ points. This is consistent with the findings in Kaplan, Kimport, Lawrence, Peilen, and Weinreich \cite{kklpw} about the coefficients of the function that counts arcs in $\mb{P}^2(\F_q)$.

We implement this algorithm in Sage \cite{sage} to count $n$-caps in $\mb{P}^3(\F_q)$ for small $n$. We arrive at the following theorem.

\begin{thm}\label{caps_thm}
	Let $$a(q) = \begin{cases}
		1 & 2 \mid q\\
		0 & 2 \nmid q.
	\end{cases} $$ Then
	\begin{align*}
		c_{3,3}(q) &=(q^2 + q + 1)(q^2 + 1)(q + 1)^2q^3\\[3pt]
		c_{4,3}(q)&=(q^3 + q^2 - 2q + 1)(q^2 + q + 1)(q^2 + 1)(q + 1)^2q^3 \\[3pt]
		c_{5,3}(q)&= (q^4 + 4q^3 + q^2 - 5q + 6)(q^2 + q + 1)(q^2 + 1)(q + 1)^2(q - 1)^2q^3\\[3pt]\end{align*}\begin{align*}
		c_{6,3}(q)&=\bigg(q^{11} + 5q^{10} - 5q^9 - 39q^8 + 35q^7 + 96q^6 - 210q^5 + 144q^4 + 169q^3 - 116q^2 \\
		&- 260q + 240\bigg)(q^2 + q + 1)(q^2 + 1)(q + 1)(q - 1)q\\[3pt]
		c_{7,3}(q) &=\bigg(q^{10} + 6q^9 - 13q^8 - 98q^7 + 148q^6 + 629q^5 - 1461q^4 - 686q^3 + 6462q^2\\
		&  - 11004q + 7470-30a(q)\bigg)(q^2 + q + 1)(q^2 + 1)(q + 1)^2(q - 1)^2q^3.
	\end{align*}
	
\end{thm}

In particular, $c_{n,3}(q)$ is a polynomial in $q$ when $n \le 6$ and is a quasipolynomial in $q$ when $n=7$. 
	

\section{Proof of Main Results}
There are no hyperfigurations on $n \le 5$ points, so the Algorithm \ref{caps_alg} produces the exact polynomial formulas for $c_{n,3}(q)$ when $n \le 5$ given in Theorem \ref{caps_thm}. 

We could count 4-caps and 5-caps combinatorially as well. For example, to count 4-caps, we need to count sets of four points so that no three are on a line. We can choose any two points, then pick a third point not on the line determined by those first two points. Lastly, we must pick a fourth point that does not lie on any of the three lines formed by the first three points. This gives the count
$$
(q^3+q^2+q+1)(q^3+q^2+q)(q^3+q^2)(q^3+q^2-2q+1)
$$
which matches the formula in Theorem \ref{caps_thm} after factoring. We can count 5-caps similarly, but the count is much more intricate with cases depending on whether or not the first four points lie on a common plane.

There is exactly one hyperfiguration on 6 points, which has planes $\{\{0,1,2,3\}, \{0,1,2,4\},$ $\{0,1,2,5\}, \{0,3,4,5\},$ $\{1,3,4,5\}, \{2,3,4,5\}\}$ and lines $\{\{0,1,2\}, \{3,4,5\}\}.$ It is important to note that these lines must be skew in $\mb{P}^3(\Fq)$. Algorithm \ref{caps_alg} gives the following expression for $c_{6,3}(q)$ in terms of a $\Z[q]$-linear combination of 1 and  $A_6(4,q)$, the number of strong realizations of this hyperfiguration on 6 points. 

\begin{prop}\label{prop_6caps}
	We have
	\begin{align*}
	c_{6,3}(q) &= q^{18} + 6q^{17} + q^{16} - 39q^{15} - 20q^{14} + 76q^{13} - 80q^{12} + 79q^{11} + 132q^{10}+ 115q^9\\
	&  - 128q^8 - 176q^7 - 133q^6 + 43q^5 + 207q^4 + 136q^3 + 20q^2 + 10A_6(4,q) - 240q.
	\end{align*}
\end{prop}

In \cite{ish_arcsp3}, we show that
$$
A_6(4,q)  = (q^2 + q + 1)(q^2 + 1)(q + 1)^2(q - 1)^2q^6.
$$
Plugging this into Proposition \ref{prop_6caps} and simplifying gives the count for $c_{6,3}(q)$ in Theorem \ref{caps_thm}.

Finally, we consider 7-caps in $\mb{P}^3(\F_q)$. There are 6 non-isomorphic hyperfigurations on 7 points, which we list below.

\begin{align*}
	h_1: \;&\mc{H} = \big\{\{0, 1, 2, 3\}, \{0, 1, 4, 5\}, \{0, 2, 4, 6\}, \{1, 2, 5, 6\}, \{1, 3, 4, 6\}, \{2, 3, 4, 5\}\big\},\\
	&\mc{L} =  \big\{\big\}\\[5pt]
	h_2: \;&\mc{H}= \big\{\{0, 1, 2, 3\}, \{0, 1, 4, 5\}, \{0, 2, 4, 6\}, \{0, 3, 5, 6\}, \{1, 2, 5, 6\}, \{1, 3, 4, 6\}, \{2, 3, 4, 5\}\big\}, \\
	&\mc{L}=\big\{\big\}\\[5pt]
	h_3:\; &\mc{H}= \big\{\{0, 1, 2, 3\}, \{0, 1, 2, 4\}, \{0, 1, 2, 5\}, \{0, 1, 2, 6\}, \{0, 3, 4, 5\}, \{1, 3, 4, 6\}, \{2, 3, 5, 6\}\big\},\\
	&\mc{L}= \big\{\{0,1,2\}\big\}\\[5pt]
	h_4:\; &\mc{H} =\big\{\{0, 1, 2, 3, 4\}, \{0, 1, 2, 5\}, \{0, 1, 2, 6\}, \{0, 3, 4, 5\}, \{0, 3, 4, 6\},\{1, 3, 5, 6\}, \{2, 4, 5, 6\}\big\} ,\\
	&\mc{L}=\big\{\{0, 1, 2\}, \{0,3, 4\}\big\}\\[5pt]
	h_5:\;&\mc{H}=\big\{\{0, 1, 2, 3, 4\}, \{0, 1, 2, 3, 5\}, \{0, 1, 2, 3, 6\}, \{0, 4, 5, 6\}, \{1, 4, 5, 6\}, \{2, 4, 5, 6\}, \{3, 4, 5,6 \}\big\},\\
	&\mc{L}=\big\{\{0, 1, 2,3\}, \{4, 5, 6\}\big\}\\[5pt]
	h_6:\; &\mc{H}= \big\{\{0,1,2,3,4,5,6\}\big\},\\
	&\mc{L}= \big\{\{0, 1, 2\}, \{0, 3, 4\}, \{0, 5, 6\}, \{1, 3, 5\}, \{1, 4, 6\}, \{2, 3, 6\}, \{2, 4, 5\}\big\}\\
\end{align*}

We run Algorithm \ref{caps_alg} and find the following expression for $c_{7,3}(q)$ in terms of a $\Z[q]$-linear combination of 1, $A_6(4,q)$, and $A_{h_i}(4,q)$ for $1\le i \le 6$. 

\begin{prop}\label{prop_7caps}
	We have
	\begin{align*}
		c_{7,3}(q) &= q^{21} + 7q^{20 }- 7q^{19} - 112q^{18} - 28q^{17} + 637q^{16} - 119q^{15}- 1840q^{14}\\
		& + 4676q^{13} - 4949q^{12} - 1436q^{11} + 4136q^{10} - 107q^9 + 7839q^8 - 5908q^7 \\
		&- 2184q^6 - 4542q^5 - 3534q^4 + 7470q^3 + (70q^3 + 70q^2- 560q+ 280) A_6(4,q)  \\
		&- 105A_{h_5}(4,q) - 30A_{h_6}(4,q).
	\end{align*}
\end{prop}

In \cite{ish_arcsp3}, we give formulas for the number of strong realizations of each hyperfiguration on seven points. There is a small error in the computation of $A_{h_6}(4,q)$, which we correct now. This error does not affect any results in \cite{ish_arcsp3}.

\begin{prop}\label{strong_real7}\cite[Proposition 3.8]{ish_arcsp3}
	Let $$a(q) = \begin{cases}1 & q \equiv 0 \pmod{2}\\0 & q \equiv 1 \pmod{2}\end{cases}.$$	The number of strong realizations for each hyperfiguration is given by
	
	\begin{align*}
		A_{h_1}(4,q)&=(1-a(q)) \cdot |\PGL_4(\F_q) |\\[3pt]
		A_{h_2}(4,q) &= a(q)  \cdot |\PGL_4(\F_q)|\\[3pt]
		A_{h_3}(3,q) &= (q-2) \cdot |\PGL_4(\F_q)|\\[3pt]
		A_{h_4}(4,q)&=  |\PGL_4(\F_q)| \\[3pt]
		A_{h_5}(4,q) &= (q^2 + q + 1)(q^2 + 1)(q + 1)^2(q - 1)^2(q - 2)q^6\\[3pt]
		A_{h_6}(4,q) &= a(q)\cdot (q^3+q^2+q+1) \cdot |\PGL_3(\F_q)|.
	\end{align*}
\end{prop}

\begin{proof}
	We give a corrected proof of the calculation for $A_{h_6}(4,q)$. Notice that all points in a strong realization of $h_6$ lie in the same plane. Fix a particular plane $H$ in $\mb{P}^3(\F_q)$. Each strong realization of $h_6$ projects to a copy of the Fano plane in $H$.  Therefore we can count strong realizations of $h_6$ in $H$ by counting strong realizations of the Fano plane in $\mb{P}^2(\F_q)$. This is well-known to be $|\PGL_3(\F_q)| \cdot a(q)$. We then multiply by the number of planes in $\mb{P}^3(\F_q)$.
\end{proof}
\begin{proof}[Proof of Theorem \ref{caps_thm}]
Apply Proposition \ref{strong_real7} to Proposition \ref{prop_7caps} and simplify.
\end{proof}


\section{Future Directions}
Our algorithm from \cite{ish_arcsp3} that counts arcs in $\mb{P}^3(\F_q)$ leads to an algorithm that counts caps in $\mb{P}^3(\F_q)$. In \cite{ish_arcsp3}, we also give a road map to count arcs in higher-dimensional projective space. It would be interesting to develop similar algorithms to count arcs and caps in $\mb{P}^4(\F_q)$. This will likely be computationally intensive as the analogue of a planar space would now contain points, lines, planes, and hyperplanes and the interactions between these objects are more complicated.

It is interesting to note that the non-polynomial behavior of $c_{7,3}(q)$ comes from counting strong realizations of $h_6$. Each strong realization of $h_6$ projects to a strong realization of the Fano plane in $\mb{P}^2(\F_q)$. Therefore we can say that the non-polynomial behavior of $c_{7,3}(q)$ comes from the Fano plane. It is unlikely that for large $n$, all hyperfigurations contributing to the non-polynomial behavior of $c_{n,3}(q)$ will project to superfigurations in the projective plane. We pose the question: what is the smallest $n$ for which a hyperfiguration exists with all points not in a common plane?

In a future paper, we will study 8-arcs and 8-caps in $\mb{P}^3(\F_q)$. It is likely that this formula is also quasipolynomial in $q$. Since the number of 10-arcs in $\mb{P}^2(\F_q)$ is non-quasipolynomial, it seems extremely likely that $c_{10,3}(q)$ is non-quasipolynomial as well. The remaining open case is whether or not $c_{9,3}(q)$ is quasipolynomial.

\section{Acknowledgments}
The author thanks Nathan Kaplan and Max Weinreich for helpful comments and suggestions.
\bibliographystyle{habbrv}
\footnotesize{\bibliography{../../Bibliography/bib_all}}

\end{document}